\documentclass{amsart}
\usepackage{amssymb}
\usepackage{amsfonts}
\usepackage{latexsym}

\newtheorem{theorem}{Theorem}[section]
\newtheorem{lemma}[theorem]{Lemma}
\newtheorem{proposition}[theorem]{Proposition}
\newtheorem{corollary}[theorem]{Corollary} 
\theoremstyle{definition}  
\newtheorem{definition}[theorem]{Definition}

\newtheorem{remark}[theorem]{Remark}

\newcommand{\Tr}{\text{Tr}}%{\text{Tr}\,} 
\newcommand{\id}{\text{id}}

\newcommand{\End}{\text{End}} 

\newcommand{\Hom}{\text{Hom}}

\newcommand{\Rep}{\text{Rep}}

\newcommand{\eps}{\varepsilon}

\newcommand{\bV}{{\bar V}}

\newcommand{\ben}{\begin{enumerate}}
\newcommand{\een}{\end{enumerate}}

\newcommand{\cC}{{\mathcal C}}
\newcommand{\cQ}{{\mathcal Q}}

\newcommand{\cV}{{\mathcal V}}

\newcommand{\cD}{{\mathcal D}}
\newcommand{\cT}{{\mathcal T}}
\newcommand{\cI}{{\mathcal I}}
\newcommand{\cM}{{\mathcal M}}
\newcommand{\be}{{\bf 1}}
\newcommand{\ba}{\bar \alpha}
\newcommand{\bb}{\bar \beta}

\newcommand{\BZ}{{\mathbb Z}}

\newcommand{\BC}{{\mathbb C}}
\newcommand{\BA}{{\mathbb A}}
\newcommand{\sq}{$\square$}

\newcommand{\Ve}{\mbox{Vec}}

\newcommand{\Fun}{\mbox{Fun}}
\newcommand{\Ima}{\mbox{Im}}
\newcommand{\Gr}{\mbox{Gr}}

\newcommand{\uloop}{\begin{picture}(10,10)
\put(0,0){\circle*{3}}
\qbezier(0,0)(-5,5)(0,10)
\qbezier(0,0)(5,5)(0,10)
\end{picture}}
\newcommand{\rloop}{\begin{picture}(10,10)
\put(0,0){\circle*{3}}
\qbezier(0,0)(5,5)(10,0)
\qbezier(0,0)(5,-5)(10,0)
\end{picture}}
\newcommand{\dloop}{\begin{picture}(10,10)
\put(0,0){\circle*{3}}
\qbezier(0,0)(-5,-5)(0,-10)
\qbezier(0,0)(5,-5)(0,-10)
\end{picture}}
\newcommand{\lloop}{\begin{picture}(10,10)
\put(0,0){\circle*{3}}
\qbezier(0,0)(-5,-5)(-10,0)
\qbezier(0,0)(-5,5)(-10,0)
\end{picture}}
\newcommand{\ulloop}{\begin{picture}(10,10)
\put(0,0){\circle*{3}}
\qbezier(0,0)(-5,0)(-10,10)
\qbezier(0,0)(-5,15)(-10,10)
\end{picture}}

\begin{document}

\title{Module categories over representations of $SL_q(2)$ and graphs}
\begin{abstract} We classify module categories over the category of
representations of quantum $SL(2)$ in a case when $q$ is not a root of
unity. In a case when $q$ is a root of unity we classify module categories
over the semisimple subquotient of the same category.
\end{abstract}

\author{Pavel Etingof}
\address{Department of Mathematics, Massachusetts Institute of Technology,
Cambridge, MA 02139, USA}
\email{etingof@math.mit.edu}

\author{Viktor Ostrik}
\address{Department of Mathematics, Massachusetts Institute of Technology,
Cambridge, MA 02139, USA}
\email{ostrik@math.mit.edu}

\maketitle

\section{Introduction}

Let $q$ be a root of unity of even order $N>4$. 
Let ${\mathcal C}_q$ be the corresponding fusion category 
of representations of the quantum group $SL_q(2)$. 
It is known (see \cite{Oc},\cite{KO} and references therein) that
indecomposable semisimple module categories over ${\mathcal C}_q$
correspond to the ADE Dynkin diagrams with Coxeter number
$h=N/2$. This fact may be viewed as the ``quantum McKay's
correspondence''. More specifically, the module categories 
in question may be viewed as ``quantum
finite subgroups in $SL_q(2)$'', by analogy with finite subgroups
of $SL(2)$, which define module categories over ${\rm
Rep}(SL(2))$ and are parametrized by the ADE affine Dynkin
diagrams by virtue of the classical McKay's correspondence.   

In this paper, we generalize this picture to the case 
of any nonzero complex number $q$, not equal to $\pm i$. 
Namely, let $q$ be such a number.
If $q=\pm 1$ or $q$ is not a root of unity, let ${\mathcal C}_q$ denote
the category of representations of the quantum group 
$SL_q(2)$. If $q$ is a root of unity such that $q^4\ne 1$, we 
let ${\mathcal C}_q$ denote the fusion category 
attached to the quantum group $SL_q(2)$. 
We classify indecomposable semisimple module
categories over ${\mathcal C}_q$ with finitely many simple
objects. It turns out that such module categories 
are parametrized by connected graphs equipped with bilinear forms 
satisfying some relations. In the case when $q$ is a root of unity of even
order $N>4$, this easily yields the classification 
of \cite{Oc},\cite{KO}; so in particular we obtain 
a very simple proof of the result of \cite{KO}, which does not
involve vertex algebras and conformal inclusions (in fact, 
this proof is close to the original approach of \cite{Oc}). 

A striking property of our classification is that while all
connected graphs do appear, trees appear only for special values
of $q$, namely such that $-q-q^{-1}$ is an eigenvalue of the
adjacency matrix that admits an eigenvector with nonvanishing
entries. Thus we discover somewhat unexpected 
``combinatorial'' peculiarities of $SL_q(2)$ at algebraic special
values of $q$ which are not roots of unity. 

On the contrary, we show 
that graphs with cycles appear for generic $q$ (or, equivalently, over
$\overline{{\mathbb C}(q)}$). This explains why 
the only finite subgroups of $SL(2)$ which admit a continuous quantum
deformation (into subgroups of $SL_q(2)$) are ${\mathbb Z}/n{\mathbb Z}$:  
for them, the corresponding affine Dynkin graph has a cycle (type
$\tilde A_{n-1}$), while for other cases (types $\tilde D_n,\tilde
E_n$), this graph is a tree. 

{\bf Acknowledgments.}
We are pleased to thank Dan Spielman for providing reference 
\cite{Fi}. The research of P.E. was partially supported by
the NSF grant DMS-9988796, and was done in part for the Clay Mathematics 
Institute. The research of V.O. was supported by the NSF 
grant DMS-0098830.  

\section{Main equation}
\subsection{Quantum $SL(2)$} We will work over the field $\BC$ of complex 
numbers. Let $q\in \BC$ be a nonzero number, $q^2\ne -1$. Recall (see e.g. 
\cite{K}) that the Hopf algebra $SL_q(2)$ 
is defined by generators $a,b,c,d$ and relations:
$$ba=qab,\; db=qbd,\; ca=qac,\; dc=qcd,\; bc=cb,$$
$$ad-da=(q^{-1}-q)bc,\; ad-q^{-1}bc=1,$$
$$\Delta \left( \begin{array}{cc}a&b\\ c&d\end{array}\right) =
\left( \begin{array}{cc}a&b\\ c&d\end{array}\right)\otimes
\left( \begin{array}{cc}a&b\\ c&d\end{array}\right),$$
$$\eps \left( \begin{array}{cc}a&b\\ c&d\end{array}\right) =
\left( \begin{array}{cc}1&0\\ 0&1\end{array}\right),\;
S \left( \begin{array}{cc}a&b\\ c&d\end{array}\right) =
\left( \begin{array}{cc}d&-qb\\ -q^{-1}c&a\end{array}\right).$$
Let $\tilde \cC_q$ denote the tensor category of finite dimensional comodules
over $SL_q(2)$. 

%Recall that since the Hopf algebra $SL_q(2)$ is cobraided 
%(see \cite{K}) the category $\tilde \cC_q$ is braided. Moreover, the category
%$\tilde \cC_q$ has a natural structure of a ribbon tensor category, see
%{\em loc. cit.} 
 
Let $\be \in \tilde \cC_q$ denote the unit object and
let $V\in \tilde \cC_q$ be a two dimensional comodule $V$ with 
the basis $x,y$ and the coaction given by 
$$\Delta_V \left( \begin{array}{c}x\\ y\end{array}\right) =
\left( \begin{array}{cc}a&b\\ c&d\end{array}\right) \otimes
\left( \begin{array}{c}x\\ y\end{array}\right) .$$
The following well known property of the object $V$ will be crucial for us:

The object $V\in \tilde \cC_q$ is selfdual, moreover for any
isomorphism $\phi :V\to V^*$ the composition
$$\be \stackrel{coev_V}{\longrightarrow}V\otimes V^*
\stackrel{\phi \otimes \phi^{-1}}{\longrightarrow}V^*\otimes V
\stackrel{ev_V}{\longrightarrow}\be \eqno(1)$$
equals to $-(q+q^{-1})\id_\be$.

Indeed, let $\delta_x, \delta_y\in V^*$ be the dual basis to $x,y\in V$. By 
the definition
$$\Delta_{V^*} \left( \begin{array}{c}\delta_x\\ \delta_y\end{array}\right) =
\left( \begin{array}{cc}d&-q^{-1}c\\ -qb&a\end{array}\right) \otimes
\left( \begin{array}{c}\delta_x\\ \delta_y\end{array}\right) .$$
It is easy to check that the map $\phi(x)=\delta_y, \phi(y)=-q\delta_x$ is
an isomorphism of comodules. Finally, the composition in (1) equals to
$$1\mapsto x\otimes \delta_x+y\otimes \delta_y\mapsto \delta_y\otimes 
(-q^{-1}y)+(-q\delta_x)\otimes x\mapsto -q^{-1}-q.$$

Note that since $V\in \tilde \cC_q$ is irreducible, the isomosphism $\phi$ is 
unique up to scaling, (so composition (1) equals $-(q+q^{-1})\id_{\be}$ 
for any $\phi$). From now on we fix a choice of such isomorphism.

Equivalently, we can replace the isomorphism $\phi$ by two maps 
$$\alpha:=(\id_V\otimes \phi^{-1})\circ coev_V:\be \to V\otimes V,\; \beta:=
ev_V\circ (\phi \otimes \id_V): V\otimes V\to \be$$ such that:

1) The compositions $V\stackrel{\alpha\otimes \id_V}{\longrightarrow}V
\otimes V\otimes V\stackrel{\id_V \otimes \beta}{\longrightarrow}V$ and 
$V\stackrel{\id_V \otimes \alpha}{\longrightarrow}V\otimes V\otimes V
\stackrel{\beta \otimes \id_V}{\longrightarrow}V$ both equal to 
$\id_V :V\to V$.

2) The composition $\be \stackrel{\alpha}{\longrightarrow}V\otimes V
\stackrel{\beta}{\longrightarrow}\be$ equals to $-(q+q^{-1})\id_\be$.

Indeed, the map $\phi$ can be reconstructed from the pair $(\alpha,\beta)$ as 
the composition $V\stackrel{\id_V\otimes coev_V}{\longrightarrow}V\otimes 
V\otimes V^*\stackrel{\beta \otimes \id_{V^*}}{\longrightarrow}V^*$.

\subsection{Turaev's construction: the generic case} Recall that in the case 
when $q$ is not a root of unity or $q=\pm 1$ the category $\tilde \cC_q$ is 
semisimple and we have a unique isomorphism of the Grothendieck rings (as 
based rings) $\Gr (\tilde \cC_q)\simeq \Gr (\tilde \cC_1)=\Gr (\Rep(SL(2))$, 
see e.g. \cite{K}. In other words, the category $\tilde \cC_q$ has exactly
one simple comodule in each dimension and tensor products of simple 
comodules are decomposed in the same way as for $SL(2)$.

In Chapter XII of \cite{T} V.~Turaev gave a topological construction of the 
category $\tilde \cC_q$. We reformulate his results in the following way:

\begin{theorem}\label{turaev1} {\em (\cite{T})} 
Assume that $q$ is not a root of unity or $q=\pm 1$.
The triple $(\tilde \cC_q,V,\phi)$ has the following universal property: let 
$\cD$ be an abelian monoidal category, let $W\in \cD$ be a 
right rigid object and $\Phi: W\to W^*$ be an isomorphism such that the
composition morphism
$$\be \stackrel{coev_W}{\longrightarrow}W\otimes W^*
\stackrel{\Phi \otimes \Phi^{-1}}{\longrightarrow}W^*\otimes W
\stackrel{ev_W}{\longrightarrow}\be \eqno(2)$$
equals to $-(q+q^{-1})\id_\be$.
Then there exists a unique tensor functor $F: \tilde \cC_q\to \cD$ such
that $F(V)=W$ and $F(\Phi )=\phi$.
\end{theorem}

{\bf Sketch of proof.} We will freely use the notation from Chapter XII of
\cite{T}. Let $\ba =(\id_W\otimes \Phi^{-1})\circ coev_W:
\be \to W\otimes W,\; \bb:=ev_W\circ (\Phi \otimes \id_W): W\otimes W\to \be$.
Obviously, the morphisms $\ba$ and $\bb$ induce the homomorphisms
$E_{k,l}\to \Hom(W^{\otimes k},W^{\otimes l})$ compatible with the compositions
(here $E_{k,l}$ are the skein modules, \cite{T} XII.1.1). In particular, for
$k=l$ we get the homomorphism $E_k\to \End(W^{\otimes k})$ where $E_k=E_{k,k}$
is the Temperley-Lieb algebra. Let $f_k\in E_k$ be the Jones-Wenzl projectors
(see \cite{T} XII.4.1). Set $a=\sqrt{-q}$ and recall that Turaev defined the
category $\cV(a)$ (see \cite{T} XII.6) objects of which are sequences
$(j_1,j_2,\ldots,j_l)\in \BZ_{\ge 0}^l$. Define the functor $\tilde F:
\cV(a)\to \cD$ by $F((j_1,j_2,\ldots,j_l))=\Ima (f_{j_1}\otimes f_{j_2}\otimes
\cdots \otimes f_{j_l})\subset W^{\otimes (j_1+j_2+\ldots +j_l)}$ and endow it
with the obvious tensor structure; since the morphisms in the category $\cV(a)$
are defined in terms of the skein modules, the functor $\tilde F$ is well
defined on morphisms. Let us apply this construction to the case $\cD =\tilde
\cC_q, W=V, \Phi =\phi$. We get the functor $\tilde F_q:\cV(a)\to \tilde
\cC_q$. The calculations in \cite{T} XII.8 show that the functor $\tilde F_q$
is an equivalence of categories. Thus we can set $F=\tilde F\circ \tilde 
F_q^{-1}$ and the Theorem is proved. \sq

\begin{remark}\label{four}
(i) One can require from $\cD$ to be only a Karoubian category,
that is an additive category where any projector has an image.

(ii) Theorem \ref{turaev1} implies immediately that the categories 
$\tilde \cC_q$ and
$\tilde \cC_{q^{-1}}$ are equivalent. Of course this is well known. Another
fact of a similar kind is the following. Let $\tilde \cC_q^{\pm}\subset
\tilde \cC_q$ be the full subcategories with objects which are direct sums
of odd/even dimensional simple comodules depending on $\pm$. Clearly,
$\tilde \cC_q=\tilde \cC_q^+\oplus \tilde \cC_q^-$. Moreover,
$\tilde \cC_q^+\otimes \tilde \cC_q^+\subset \tilde \cC_q^+$,
$\tilde \cC_q^+\otimes \tilde \cC_q^-\subset \tilde \cC_q^-$,
$\tilde \cC_q^-\otimes \tilde \cC_q^+\subset \tilde \cC_q^-$,
$\tilde \cC_q^-\otimes \tilde \cC_q^-\subset \tilde \cC_q^+$.
In other words, the category $\tilde \cC_q$ is $\BZ/2\BZ-$graded.  
Now one can twist the associativity isomorphism in $\tilde \cC_q$ by
changing the sign of the associativity isomorphism $(X\otimes Y)\otimes Z\to
X\otimes (Y\otimes Z)$ when $X,Y,Z \in \tilde \cC_q^-$. Let us denote the
category twisted in such a way by $\tilde \cC_q^{tw}$. It follows from
Theorem \ref{turaev1} that the category $\tilde \cC_q^{tw}$ is equivalent to
$\tilde \cC_{-q}$. This is also well known; moreover both facts above
remain true when $q$ is a root of unity.

(iii) One should be very careful with universal properties of tensor 
categories: for example the universal category with an object $X$ such that
$X\otimes X=\be$ does not exist. We expect that when $q$ is a root of unity
the universal abelian category in a sense of Theorem \ref{turaev1} does not
exist. In contrast the universal Karoubian category clearly exists 
(it coincides with Karoubian envelope of the skein category from \cite{T} 
XII.2) and coincides with the category of tilting modules $\cT_q\subset 
\tilde \cC_q$ (this is a consequence of the quantized Schur-Weyl duality, see
\cite{D}).

(iv) One can restate Theorem \ref{turaev1} in the following way: the tensor
functors $F:\tilde \cC_q\to \cD$ are in the one to one correspondence with
objects $W\in \cD$ together with an isomorphism $\Phi: W\to W^*$ such that
the composition (2) equals to $-q-q^{-1}$.
\end{remark}

\subsection{Turaev's construction: the roots of unity case}
In the case when $q$ is a root of unity of order $N\ge 3$ the category 
$\tilde \cC_q$ is not semisimple, see \cite{K}. Let $\cT_q$ denote the full 
additive (nonabelian) subcategory of $\tilde \cC_q$ whose objects are direct 
summands of $V^{\otimes n}$; clearly the subcategory $\cT_q$ is closed under 
the tensor product (the category $\cT_q$ is the category of tilting modules,
see e.g. \cite{BK}). It is well known that the additive subcategory $\cI_q$ 
of $\cT_q$ generated by indecomposable modules of zero quantum dimension
is a tensor ideal and thus the quotient $\cC_q=\cT_q/\cI_q$ is a well defined
semisimple tensor category, see \cite{BK}. The object $V\in \tilde \cC_q$ can
be considered as an object of $\cC_q$ and for the isomorphism $\phi : V\to V^*$
the composition (2) equals to $-q-q^{-1}$. But the universal property of the
category $\cC_q$ is a little bit more delicate. Let $\cD$ be an abelian
monoidal category, $W\in \cD$ be a right rigid object and $\Phi :W\to W^*$
be an isomorphism such that the composition (2) equals to $-q-q^{-1}$.
In the same way as in the discussion of Theorem \ref{turaev1} we have
homomorphisms $E_{k,l}\to \Hom(W^{\otimes k}, W^{\otimes l})$ where $E_{k,l}$
are the skein modules. Set $N_*=N$ if $N$ is odd and $N_*=N/2$ if $N$ is even.
Recall that the last Jones-Wenzl idempotent which is possible to define is
$f_{N_*-1}$, see \cite{T} XII.4.3.

\begin{theorem}\label{turaev2} Assume that $q$ is a root of unity of order
$N\ge 3$. The triple $(\cC_q,V,\phi)$ has the following universal property:
let $\cD$ be an abelian monoidal category, let $W\in \cD$ be a 
right rigid object and $\Phi: W\to W^*$ be an isomorphism such that the
composition morphism
$$\be \stackrel{coev_W}{\longrightarrow}W\otimes W^*
\stackrel{\Phi \otimes \Phi^{-1}}{\longrightarrow}W^*\otimes W
\stackrel{ev_W}{\longrightarrow}\be \eqno(2)$$
equals to $-(q+q^{-1})\id_\be$. In addition let us assume that the image of
$f_{N_*-1}$ in $\Hom(W^{\otimes {N_*-1}}, W^{\otimes {N_*-1}})$ is zero.
Then there exists a unique tensor functor $F: \cC_q\to \cD$ such
that $F(V)=W$ and $F(\Phi )=\phi$.

In other words, the tensor functors $F: \cC_q\to \cD$ are in bijection with
objects $W\in \cD$ together with an isomorphism $\Phi :W\to W^*$ such that
the composition (2) equals to $-q-q^{-1}$ and the image of
$f_{N_*-1}$ in $\Hom(W^{\otimes {N_*-1}}, W^{\otimes {N_*-1}})$ is zero.
\end{theorem}

{\bf Remarks on proof.} As it is explained in Remark \ref{four} (iii) we 
have a unique tensor functor $F: \cT_q\to \cD$ such that $F(V)=W$ and
$F(\Phi)=\phi$. So one just needs to check that this functor maps $\cI_q$ to
zero. But it is well known (see e.g. \cite{BK}) that the tensor ideal $\cI_q$ 
is generated by $\Ima (f_{N_*-1})$ in a sense that any object of $\cI_q$ is
isomorphic to a direct summand of $\Ima (f_{N_*-1})\otimes T$ where $T\in 
\cT_q$. \sq

\begin{remark} In \cite{T} only the case of even $N$ is considered. We note
that the construction of \cite{T} works without any change for odd $N$ as
well. The only difference is that the resulting semisimple category is not
modular.
\end{remark}

\subsection{Main equation} Define
$$\cC_q=\left\{ \begin{array}{ll}\tilde \cC_q&\mbox{if}\; q \; \mbox{is not a
root of unity or}\; q=\pm 1;\\ 
\cT_q/\cI_q&\mbox{if}\; q \; \mbox{is a root of unity},\; q\ne \pm 1.
\end{array}\right.$$
The aim of this note is to classify the semisimple module categories with 
finitely many simple objects over the category $\cC_q$. Here is our main
result:

\begin{theorem}\label{main} {\em (i)} Assume that $q$ is not a root of unity
or $q=\pm 1$. The semisimple module categories with finitely many simple 
objects over the category $\cC_q$ are classified by the following data:

1) A finite set $I$;

2) A collection of finite dimensional vector spaces $V_{ij}$, $i,j\in I$;

3) A collection of nondegenerate bilinear forms $E_{ij}: V_{ij}\otimes V_{ji}
\to \BC$,

\noindent subject to the following condition: for each $i\in I$ we have
$$\sum_{j\in I}\Tr (E_{ij}(E_{ji}^t)^{-1})=-q-q^{-1}.\eqno(3)$$
\end{theorem}

\begin{proof} Let $\cM$ be a semisimple module category over $\cC_q$ with
finitely many simple objects. Let $I$ be the set
of the isomorphism classes of simple objects in $\cM$. The structure of module
category on $\cM$ is the same as the tensor functor 
$F:\cC_q\to \Fun(\cM,\cM)$ where $\Fun(\cM,\cM)$ is the category of additive
functors from $\cM$ to itself, see \cite{O}. Recall that the category 
$\Fun(\cM,\cM)$ is
identified with the category of $I\times I-$graded vector spaces with obvious
``matrix'' tensor product. By Remark \ref{four} (iv) the functors $F:\cC_q\to
\Fun(\cM,\cM)$ are bijective to the objects $\bV =(V_{ij})\in \Fun(\cM,\cM)$
together with an isomorphism $\Phi: \bV \to \bV^*=(V_{ji}^*)$ (equivalently,
this is a collection of nondegenerate bilinear forms $E_{ij}: V_{ij}\otimes 
V_{ji}\to \BC$) such that the morphism (2) equals to $-q-q^{-1}$. It is
obvious that the last condition is equivalent to the condition (3). The
theorem is proved. \end{proof}

\section{Solutions of the main equation} 

\subsection{Deformations} Let us fix a finite set $I$, the
numbers $a_{ij}=\dim V_{ij}=\dim V_{ji}$ and try to analyze the corresponding
solutions of the main equation. It is convenient to represent these data
as a graph $\Gamma =(I, \{ a_{ij}\} )$ with the set of vertices $I$ and 
$a_{ij}$ edges joining the vertices $i$ and $j$ (the matrix $A=(a_{ij})$ is 
the adjacency matrix of this graph). We are going to classify the graphs 
with respect to the deformation behavior of solutions of the main equation.
If we fix the vector spaces $V_{ij}$ of dimensions $a_{ij}$, the set of 
solutions of equation (3) is clearly an affine algebraic variety ${\mathfrak
M}$. The group $G=\prod_{i,j}GL(V_{ij})$ acts naturally on ${\mathfrak M}$.
Since we are interested in the solutions of the main equation only up to
isomorphism we define the set of solutions of the main equation to be the
set ${\mathfrak M}/G$ of orbits of $G$ on ${\mathfrak M}$. In general
${\mathfrak M}/G$ has no structure of an algebraic variety; so let
${\mathfrak M}//G$ denote the quotient in the sense of the invariant theory,
that is ${\mathfrak M}//G$ is the set of closed $G-$orbits on ${\mathfrak M}$.
Now ${\mathfrak M}//G$ has a structure of an algebraic variety; we will see
that the natural map ${\mathfrak M}/G \to {\mathfrak M}//G$ is finite to one
and is one to one on an open nonempty subset of ${\mathfrak M}//G$. Thus we
define the dimension of the set of solutions of the main equation to be equal
to the dimension of the variety ${\mathfrak M}//G$. In such situation we will
say that ${\mathfrak M}/G$ is a moduli space (even if it is not an algebraic
variety). 

\begin{definition} (i) We will say that a graph is super-rigid if the 
main equation (3) admits only finitely many solutions for finitely many 
values of $q$ and no solutions for other values of $q$.

(ii) We will say that a graph is rigid if the main equation (3) admits
only finitely many solutions for all but finitely many values of $q$.

(iii) We will say that a graph is non-rigid if it is not rigid.
\end{definition}

\begin{remark} One says that a graph is strictly rigid if the main equation
admits only finitely many solutions for all values of $q$. We will see later
that the graph 
\begin{picture}(20,10)
\put(0,0){\circle*{3}}
\put(0,0){\line(1,0){20}}
\put(20,0){\circle*{3}}
\put(0,0){\line(1,1){10}}
\put(10,10){\circle*{3}}
\put(10,10){\line(1,-1){10}}
\end{picture}
is rigid. On the other hand it is easy to see that it is not strictly rigid:
for $q+q^{-1}=1$ it admits infinitely many solutions of the main
equation. \end{remark}

If a graph $\Gamma$ is a disjoint union of two subgraphs, the
cooresponding module category over $\cC_q$ is clearly a direct sum of module
subcategories corresponding to the subgraphs. Thus from now on we will study 
only connected graphs. For a graph $\Gamma$ we define its 
{\em underlying simply laced graph} as a graph $\bar \Gamma$ with the same set
of vertices and the vertices $i\ne j$ are joined by exactly one edge if 
$v_{ij}\ne 0$ and are not joined otherwise (in particular $\bar \Gamma$ has
no self-loops).

 For a graph $\Gamma =(I,\{ a_{ij}\} )$ we define the generalized number of
cycles $L(\Gamma)$ by formula
$$L(\Gamma)=\frac{1}{2}\sum_{i\ne j}a_{ij}+\sum_i[\frac{a_{ii}}{2}]-|I|+1$$
where $[\cdot]$ denotes the integer part. Note that in a case when $\Gamma$
has no self-loops $L(\Gamma)$ is just the number of loops in $\Gamma$.
We will see later that the number $L(\Gamma)$ is the expected dimension (that 
is, the difference of the number of variables and the number of equations) of 
the set of solutions of the main equation. Moreover, we will see that 
$L(\Gamma)$ coincides with (properly understood) dimension of the set of
solutions of the main equation. 

\begin{definition} (i) A connected graph is called a generalized tree
if $L(\Gamma)=0$.

(ii) A connected graph is called a 1-loop graph if $L(\Gamma)=1$.
\end{definition}

\begin{remark} (i) A connected graph is a generalized tree 
if $a_{ij}\le 1$ and its underlying simply laced graph is a tree (note that the
possibility $a_{ii}\ne 0$ is allowed).

(ii) A connected graph is a 1-loop graph if either 
$a_{ij}\le 1$ and its underlying simply laced graph has exactly $|I|$ edges or
its underlying simply laced graph is a tree, $a_{ij}\le 3$, $a_{ij}\ge 2$ for 
exactly one pair of vertices $i,j$ and $a_{ij}=3$ implies $i=j$.
\end{remark}

\begin{theorem}\label{rigid}
{\em (i)} For any graph there exists a solution of the main
equation with some $q\ne \pm i$.

{\em (ii)} A connected graph is super-rigid iff it is a generalized tree.

{\em (iii)} A connected graph is rigid but not super-rigid iff it is a 
1-loop graph.
\end{theorem}

\begin{proof}

1) A quadruple $(V_{ij},V_{ji},E_{ij},E_{ji})$ consisting of two vector
spaces $V_{ij}, V_{ji}$ and two nondegenerate bilinear forms $E_{ij}:
V_{ij}\otimes V_{ji}\to \BC$ and $E_{ji}: V_{ji}\otimes V_{ij}\to \BC$ is
isomorphic to the quadruple $(V_{ij},V_{ij}^*,\langle \cdot ,\cdot \rangle ,
\langle S\cdot ,\cdot \rangle )$ where $\langle \cdot ,\cdot \rangle :V_{ij}
\otimes V_{ij}^*\to \BC$ is the canonical pairing and $S: V_{ij}\to V_{ij}$
is an invertible linear operator; two such quadruples are isomorphic if and
only if the corresponding operators $S$ have the same Jordan form. Thus
the moduli space $\cQ(a_{ij})$ of such quadruples with $\dim (V_{ij})=\dim
(V_{ji})=a_{ij}$ has dimension 
%\footnote{We mean here the following: 
%there is a natural map from the set of equivalence classes of quadruples to 
%$\cQ(a_{ij})=\BC^*\times \ldots \times \BC^*/S_{a_{ij}}$ ($a_{ij}$ factors,
%$S_{a_{ij}}$ is the symmetric group acting by the permutations of factors) 
%which associates to a quadruple the eigenvalues
%of the operator $S$; this map is finite to one and is one to one on the
%nonempty open subset of $\cQ(a_{ij})$; main equation (3) involves only
%quantities which are regular functions on $\cQ(a_{ij})$.} 
$a_{ij}$.
The image of the map $\cQ(a_{ij})\to \BA^2$, $x=\Tr (E_{ij}(E_{ji}^t)^{-1}),
y=\Tr (E_{ji}(E_{ij}^t)^{-1})$ depends on $a_{ij}$: if $a_{ij}=1$ this is
the hyperbola $xy=1$; if $a_{ij}=2$ this is $(\BA^2-\{ xy=0\} )\cup (0,0)$;
if $a_{ij}\ge 3$ this is $\BA^2$.

2) Recall here the classification of nondegenerate bilinear forms, 
see \cite{B}, 5.6: any
pair $(V,E)$ consisting of a vector space $V$ and nondegenerate bilinear
form $E: V\otimes V\to \BC$ is up to isomorphism uniquely determined by 
the operator $S_E=E(E^t)^{-1}: V^*\to V^*$; the operator $S_E$ is conjugated
to $S_E^{-1}$ and moreover the number of Jordan cells of size $k$ with
eigenvalue $(-1)^k$ is even. Thus the moduli space $\tilde \cQ(a)$ of such 
pairs with $\dim (V)=a$ has dimension $[a/2]$ (in the same sense as before); 
the image of the map $\tilde \cQ(a)\to \BA^1$, $(V,E)\mapsto 
\Tr (E(E^t)^{-1})$ is just the point $1$ for $a=1$ and the entire $\BA^1$ for
$a\ge 2$.  

Thus we see that $L(\Gamma)$ is really the expected dimension of the set of
solutions of the main equation. It is clear that the actual dimension of the
set of solutions of the main equation is greater or equal to
$L(\Gamma)$ if this set is nonempty..

Now let us show that for any choice of the graph $\Gamma$ there 
exists a solution of the main equation. Let $(r_i)_{i\in I}$ be an 
eigenvector of the matrix $A=(a_{ij})$ with eigenvalue $\lambda$ and such that 
$\prod_{i\in I}r_i\ne 0$ (such eigenvector exists, for example one can take 
the Frobenius-Perron eigenvector). Now choose bilinear forms $E_{ij}$ in such
a way that $\Tr (E_{ij}(E_{ji}^t)^{-1})=a_{ij}r_j/r_i$ (this is possible in 
view of the remarks above). It is clear that in this way we get a solution
of the main equation with $\lambda =-q-q^{-1}$. If $\lambda$ is the
Frobenius-Perron eigenvalue we have $\lambda >0$ and thus $q\ne \pm i$.
Thus (i) is proved.

Now assume that the graph $\Gamma$ is not a generalized tree. Thus either 
$a_{i_0j_0}\ge
2$ for some $i_0, j_0\in I$, or our graph contains a cycle of length $M\ge 3$.

Case 1: $i_0=j_0$. Consider the matrix $\tilde A(u)=(\tilde a_{ij})$ where
$\tilde a_{ij}=a_{ij}$ except $\tilde a_{i_0i_0}=u\in \BC$. For real positive 
$u$ the Frobenius-Perron eigenvalue of the matrix $\tilde A(u)$ depends
nontrivially on $u$ since $\Tr (\tilde A(u))=u+\Tr(A)-
a_{i_0i_0}$. Thus for generic $u$ the matrix $\tilde A(u)$ has an eigenvector 
$(r_i(u))_{i\in I}$ with $\prod_{i\in I}r_i\ne 0$ and with an eigenvalue 
$\lambda (u)$ depending nontrivially on $u$. Thus $\lambda (u)$ takes all
values from $\BC$ except finitely many. Now a choice of $E_{ij}$ such that 
$\Tr (E_{ij}(E_{ji}^t)^{-1})=\tilde a_{ij}(u)r_j(u)/r_i(u)$ (this choice is
possible by 1) and 2) above) gives a solution of the main equation with 
$-q-q^{-1}=\lambda (u)$. Thus our graph is not super-rigid.

Case 2: $i_0\ne j_0$. In this case consider the matrix $\tilde A(u)=(\tilde 
a_{ij})$ where $\tilde a_{ij}=a_{ij}$ except $\tilde a_{i_0j_0}=u\in \BC$. 
Since $\Tr (\tilde A(u)^2)$ depends on $u$ nontrivially the same arguments
as above show that our graph is not super-rigid.

Case 3: the graph has a cycle of length $M\ge 3$. Let $(i_0,j_0)$ be an edge
from the cycle. Consider the matrix $\tilde A(u)=(\tilde 
A_{ij})$ where $\tilde A_{ij}=A_{ij}$ except $\tilde a_{i_0j_0}=u\in \BC$ and
$\tilde a_{j_0i_0}=u^{-1}$. Now $\Tr (\tilde A(u)^M)$ depends on $u$ 
nontrivially and our graph is not super-rigid.

Now we are going to prove that a generalized tree is super-rigid.

\begin{definition} An eigenvalue $\lambda$ of the matrix $A=(a_{ij})$ is 
called nondegenerate if there exists a $\lambda-$eigenvector $(r_i)_{i\in I}$ 
such that $\prod_{i\in I}r_i\ne 0$. \end{definition}

The following lemma is well known in graph theory, see \cite{Fi}. We give a
proof for the reader's convenience.

\begin{lemma}\label{fie}
(i) For any matrix $A$ with nonnegative integer entries there
exists a nondegenerate eigenvalue.

(ii) If $A$ corresponds to a 
generalized tree then a nondegenerate eigenvalue has multiplicity 1. 
\end{lemma}

\begin{proof} (i) The Frobenius-Perron eigenvalue (and its Galois conjugates) 
is always nondegenerate.

(ii) Let $\lambda$ be an eigenvalue of the matrix $A$. We are going to prove
that an $\lambda-$eigenvector $(r_i)_{i\in I}$ with $\prod_{i\in I}r_i\ne 0$
is unique up to proportionality if it exists. This would imply the statement
of Lemma since a small perturbation preserves the property 
$\prod_{i\in I}r_i\ne 0$. 

The vector $(r_i)_{i\in I}$ satisfies 
$$\sum_{i\ne j,a_{ij}=1}r_i=\left\{ \begin{array}{cc}\lambda r_j&\mbox{if}\; 
a_{jj}=0,\\
(\lambda -1)r_j&\mbox{if}\; a_{jj}=1,\end{array}\right.$$ 
(the sum is over all edges of the underlying simply laced graph with 
vertex $j$).
Let us introduce new variables parametrized by the oriented edges of the
underlying simply laced graph, $y_{ij}=r_i/r_j$. These variables satisfy   
$$y_{ij}y_{ji}=1,\; \; \;
\sum_{i\ne j,a_{ij}=1}y_{ij}=\left\{ \begin{array}{cc}\lambda &\mbox{if}\; 
a_{jj}=0,\\
\lambda -1&\mbox{if}\; a_{jj}=1.\end{array}\right. \eqno(4)$$ 
Now the result is a consequence of the following

{\bf Sublemma.} For any choice of $(\lambda_j)_{j\in I}$ the system of
equations
$$y_{ij}y_{ji}=1,\; \; \sum_{i\ne j,a_{ij}=1}y_{ij}=\lambda_j \eqno (5)$$
has at most one solution.

{\bf Proof of Sublemma.} The proof is by induction in $|I|$.
Choose a vertex $j_0$ of valency 1 of the underlying 
simply laced graph (which is a tree). Then there is only one variable 
$y_{ij_0}$ and it is uniquely defined from the equation 
$y_{ij_0}=\sum_{i\ne j_0,a_{ij_0}=1}y_{ij_0}=\lambda_{j_0}.$ 
Then the variable $y_{j_0i}=1/y_{ij_0}$ is also uniquely defined and all
other variables satisfy the system of equations of the same form with 
smaller $|I|$. The Sublemma and the Lemma are proved.
\end{proof}

Observe that in the case of a generalized tree the main equation has exactly
the form of system (4). Thus it is obvious that the only possible values of
$\lambda$ are the nondegenerate eigenvalues of $A$. So the Sublemma implies 
that a generalized tree is super-rigid. Thus (ii) is proved. 

Now we claim that for any graph $\Gamma$ the dimension of the space of 
solutions of the main equation is exactly $L(\Gamma)$. 
%(this set is not an
%algebraic variety, but there is a finite to one map from this set to an
%algebraic variety which is one to one on an open set). 
Indeed, let us
choose a spanning tree $T$ of the underlying simply laced graph $\bar \Gamma$.
Now let us choose any values of parameters attached to all edges not belonging
to $T$; in particular for any edge $ij$ from $T$ choose any values of 
$a_{ij}-1$ eigenvalues of the matrix $S$ (see 1) above). Thus we have chosen
$L(\Gamma)$ parameters. Now the main equation reduces to the system of the
shape (5) for the rest of parameters (we have one parameter for each edge
of the tree $T$). Now the Sublemma implies that we have only finitely
many solutions for these parameters. Henceforth we see that the expected
dimension $L(\Gamma)$ coincides with the actual dimension (understood as it
is explained above) of the set of solutions of the main equation. 

Now it is clear that if a graph is a 1-loop graph if and only if it is rigid
(indeed, the set of solutions of the main equation has dimension 1 and it
maps dominantly under the projection to the variable $q$).
The Theorem is proved.
\end{proof}

\begin{corollary} For any value of $|I|$ there are only finitely many rigid
(and hence super-rigid) graphs.\end{corollary}

Recall that the ultraspherical polynomials $P_n(x)$ are defined recursively
by
$$P_1(x)=1,\; P_2(x)=x,\; P_{n+2}(x)=xP_{n+1}(x)-P_n(x),\; n\ge 1.$$
It is a classical fact that the classes of simple objects in $\Gr (\cC_q)$
are given by $P_n([V])$, see e.g. \cite{K}.

\begin{corollary} Let $A$ be an indecomposable symmetric matrix with 
nonnegative entries. Then either $P_n(A)=0$ for some $n$ (all such matrices
are explicitly known and are classified by ADET graphs, see below) or $P_n(A)$
has nonnegative entries for all $n$.
\end{corollary}

\begin{proof} Let $\lambda$ be the Frobenius-Perron eigenvalue of $A$. All
indecomposable matrices with $\lambda <2$ are classified and it is well known
that for such matrices $P_n(A)=0$ for some $n$, see e.~g. \cite{EK}. Hence we 
can assume that $\lambda \ge 2$. Then the construction from the proof of 
Theorem
\ref{rigid} gives us a module category over $\cC_q$ with $q+q^{-1}=\lambda$
(and thus $q$ is not a root of unity) where the object $V$ is represented by 
the vector space valued matrix of dimension $A$. Since $P_n(A)$ gives an
action of some object in $\cC_q$, it is nonnegative. \end{proof}

Consider the case $q=1$. In this case $\cC_1\cong \Rep (SL(2))$. Any finite
subgroup $G\subset SL(2)$ gives rise to a module category $\Rep (G)$ over
$\Rep (SL(2))$. As it is known from the McKay correspondence, the corresponding
graph is then an affine ADE Dynkin diagram. Observe that the graphs of type
$\tilde D_n, \tilde E_n$ are super-rigid while the graph of type $\tilde A_n$
is just rigid. This explains the fact that among finite subgroups of $SL(2)$ 
only the cyclic subgroups corresponding to $\tilde A_n$ admit a continuos 
deformation in the ``quantum'' direction.

\subsection{Examples} In this section we will assume that $q$ is not a root
of unity. Recall that the module categories over $\cC$ with one 
simple object are the same as fiber functors ($=$ tensor functors $\cC \to 
\Ve$). We see from Theorem \ref{main} that the fiber functors on $\cC_q$
are classified by a vector space $V$ and a bilinear form $E: V\otimes V\to \BC$
such that $\Tr (E(E^t)^{-1})=-q-q^{-1}$. This is exactly the result of
J.~Bichon, see \cite{Bi} who classified all Hopf algebras $H$ such that the
category of comodules over $H$ is tensor equivalent to $\cC_q$. Thus our
Theorem \ref{main} can be considered as a generalization of Bichon's result:
we classify all weak Hopf algebras $H$ such that the category of comodules
over $H$ is tensor equivalent to $\cC_q$.
 
Observe that a graph with one vertex is rigid
iff $\dim (V)\le 3$ and is super-rigid iff $\dim (V)=1$ (the last case gives 
$q$ which is a primitive root of unity of order 3).
Here is a list of rigid graphs with $|I|=2$ (we wrote possible values of 
$q+q^{-1}$ over the super-rigid graphs):

\begin{picture}(25,10)
\put(0,0){\uloop}
\put(0,0){\dloop}
\put(0,0){\lloop}
\put(0,0){\line(1,0){10}}
\put(10,0){\circle*{3}}
\end{picture}
\begin{picture}(25,10)
\put(0,0){\uloop}
\put(0,0){\dloop}
\put(0,0){\lloop}
\put(0,0){\line(1,0){10}}
\put(10,0){\circle*{3}}
\put(10,0){\rloop}
\end{picture}
\begin{picture}(20,10)
\put(0,0){\uloop}
\put(0,0){\dloop}
\put(0,0){\line(1,0){10}}
\put(10,0){\circle*{3}}
\end{picture}
\begin{picture}(25,10)
\put(0,0){\uloop}
\put(0,0){\dloop}
\put(0,0){\line(1,0){10}}
\put(10,0){\circle*{3}}
\put(10,0){\rloop}
\end{picture}
\begin{picture}(15,10)
\put(0,0){\rloop}
\put(10,0){\circle*{3}}
\end{picture}
\begin{picture}(30,10)
\put(0,0){\rloop}
\put(10,0){\circle*{3}}
\put(10,0){\rloop}
\end{picture}
\begin{picture}(25,10)
\put(0,0){\rloop}
\put(0,0){\lloop}
\put(10,0){\circle*{3}}
\put(10,0){\rloop}
\end{picture}
\begin{picture}(35,10)
\put(0,0){\circle*{3}}
\put(0,0){\line(1,0){10}}
\put(10,0){\circle*{3}}
\put(0,3){$\pm 1$}
\end{picture}
\begin{picture}(40,10)
\put(0,0){\dloop}
\put(0,0){\circle*{3}}
\put(0,0){\line(1,0){10}}
\put(10,0){\circle*{3}}
\put(-20,3){$(-1\pm \sqrt{5})/2$}
\end{picture}
\begin{picture}(30,10)
\put(0,0){\dloop}
\put(10,0){\dloop}
\put(0,0){\circle*{3}}
\put(0,0){\line(1,0){10}}
\put(10,0){\circle*{3}}
\put(0,3){$-2$}
\end{picture}

\bigskip

Now we are going to discuss the case of super-rigid graphs. Let $\Gamma =(I,
\{ a_{ij}\} )$ be a generalized tree.
We have the following consequence of Lemma 
\ref{fie}:

\begin{proposition}\label{nond}
A solution of the main equation (3) for a generalized
tree exists if and only if $-q-q^{-1}$ is a nondegenerate eigenvalue of
$A=(a_{ij})$. In such a case the solution is unique.\end{proposition} 

Observe that since an eigenvalue of symmetric matrix is real we have

\begin{corollary} If the category $\cC_q$ has a module category corresponding
to a super-rigid graph then either $|q|=1$ or $q$ is real.
\end{corollary}  

We present here a few examples. Here is a list of generalized trees with
$\le 4$ vertices; under any graph we wrote possible values of $\lambda 
=-q-q^{-1}$ or an algebraic equation for $\lambda$; we omitted all values of 
$q$ being a root of unity (thus some graphs are omitted too; see the next 
subsection for them).

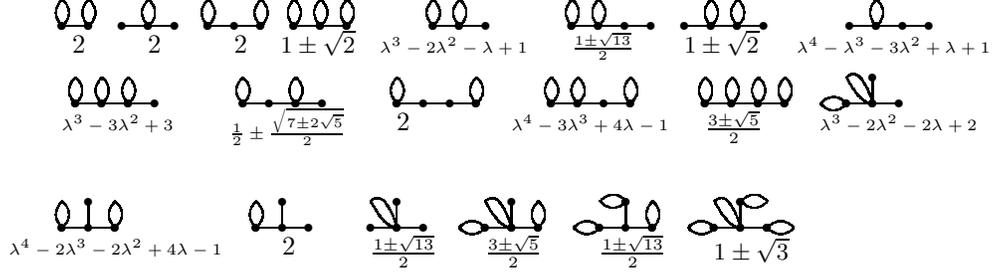
\begin{picture}(20,10)
\put(0,0){\circle*{3}}
\put(0,0){\line(1,0){10}}
\put(10,0){\circle*{3}}
\put(0,0){\uloop}
\put(10,0){\uloop}
\put(4,-10){$2$}
\end{picture}
\begin{picture}(30,10)
\put(0,0){\circle*{3}}
\put(0,0){\line(1,0){10}}
\put(10,0){\circle*{3}}
\put(10,0){\uloop}
\put(10,0){\line(1,0){10}}
\put(20,0){\circle*{3}}
\put(10,-10){$2$}
\end{picture}
\begin{picture}(30,10)
\put(0,0){\circle*{3}}
\put(0,0){\line(1,0){10}}
\put(10,0){\circle*{3}}
\put(0,0){\uloop}
\put(20,0){\uloop}
\put(10,0){\line(1,0){10}}
\put(20,0){\circle*{3}}
\put(10,-10){$2$}
\end{picture}
\begin{picture}(30,10)
\put(0,0){\circle*{3}}
\put(0,0){\line(1,0){10}}
\put(10,0){\circle*{3}}
\put(0,0){\uloop}
\put(10,0){\uloop}
\put(20,0){\uloop}
\put(10,0){\line(1,0){10}}
\put(20,0){\circle*{3}}
\put(-5,-10){\small $1\pm \sqrt{2}$}
\end{picture}
\begin{picture}(70,10)
\put(20,0){\begin{picture}(30,10)
\put(0,0){\circle*{3}}
\put(0,0){\line(1,0){10}}
\put(10,0){\circle*{3}}
\put(0,0){\uloop}
\put(10,0){\uloop}
\put(10,0){\line(1,0){10}}
\put(20,0){\circle*{3}}\end{picture}}
\put(0,-10){\tiny $\lambda^3-2\lambda^2-\lambda+1$}
\end{picture}
\begin{picture}(40,10)
\put(0,0){\circle*{3}}
\put(0,0){\line(1,0){10}}
\put(10,0){\circle*{3}}
\put(0,0){\uloop}
\put(10,0){\uloop}
\put(10,0){\line(1,0){10}}
\put(20,0){\circle*{3}}
\put(20,0){\line(1,0){10}}
\put(30,0){\circle*{3}}
\put(0,-10){\tiny $\frac{1\pm \sqrt{13}}{2}$}
\end{picture}
\begin{picture}(40,10)
\put(0,0){\circle*{3}}
\put(0,0){\line(1,0){10}}
\put(10,0){\circle*{3}}
\put(20,0){\uloop}
\put(10,0){\uloop}
\put(10,0){\line(1,0){10}}
\put(20,0){\circle*{3}}
\put(20,0){\line(1,0){10}}
\put(30,0){\circle*{3}}
\put(0,-10){\small $1\pm \sqrt{2}$}
\end{picture}
\begin{picture}(70,10)
\put(20,0){\begin{picture}(30,10)
\put(0,0){\circle*{3}}
\put(0,0){\line(1,0){10}}
\put(10,0){\circle*{3}}
\put(10,0){\uloop}
\put(10,0){\line(1,0){10}}
\put(20,0){\circle*{3}}
\put(20,0){\line(1,0){10}}
\put(30,0){\circle*{3}}
\end{picture}}
\put(0,-10){\tiny $\lambda^4-\lambda^3-3\lambda^2+\lambda +1$}
\end{picture}

\bigskip

\begin{picture}(60,10)
\put(5,0){\begin{picture}(30,10)
\put(0,0){\circle*{3}}
\put(0,0){\line(1,0){10}}
\put(10,0){\circle*{3}}
\put(0,0){\uloop}
\put(10,0){\uloop}
\put(20,0){\uloop}
\put(10,0){\line(1,0){10}}
\put(20,0){\circle*{3}}
\put(20,0){\line(1,0){10}}
\put(30,0){\circle*{3}}
\end{picture}}
\put(0,-10){\tiny $\lambda^3-3\lambda^2+3$}
\end{picture}
\begin{picture}(60,10)
\put(5,0){\begin{picture}(30,10)
\put(0,0){\circle*{3}}
\put(0,0){\line(1,0){10}}
\put(10,0){\circle*{3}}
\put(0,0){\uloop}
\put(20,0){\uloop}
\put(10,0){\line(1,0){10}}
\put(20,0){\circle*{3}}
\put(20,0){\line(1,0){10}}
\put(30,0){\circle*{3}}
\end{picture}}
\put(0,-13){\tiny $\frac{1}{2}\pm \frac{\sqrt{7\pm 2\sqrt{5}}}{2}$}
\end{picture}
\begin{picture}(40,10)
\put(0,0){\circle*{3}}
\put(0,0){\line(1,0){10}}
\put(10,0){\circle*{3}}
\put(0,0){\uloop}
\put(30,0){\uloop}
\put(10,0){\line(1,0){10}}
\put(20,0){\circle*{3}}
\put(20,0){\line(1,0){10}}
\put(30,0){\circle*{3}}
\put(0,-10){$2$}
\end{picture}
\begin{picture}(70,10)
\put(15,0){\begin{picture}(30,10)
\put(0,0){\circle*{3}}
\put(0,0){\line(1,0){10}}
\put(10,0){\circle*{3}}
\put(0,0){\uloop}
\put(10,0){\uloop}
\put(30,0){\uloop}
\put(10,0){\line(1,0){10}}
\put(20,0){\circle*{3}}
\put(20,0){\line(1,0){10}}
\put(30,0){\circle*{3}}
\end{picture}}
\put(0,-10){\tiny $\lambda^4-3\lambda^3+4\lambda -1$}
\end{picture}
\begin{picture}(50,10)
\put(0,0){\circle*{3}}
\put(0,0){\line(1,0){10}}
\put(10,0){\circle*{3}}
\put(0,0){\uloop}
\put(10,0){\uloop}
\put(20,0){\uloop}
\put(30,0){\uloop}
\put(10,0){\line(1,0){10}}
\put(20,0){\circle*{3}}
\put(20,0){\line(1,0){10}}
\put(30,0){\circle*{3}}
\put(0,-12){\small $\frac{3\pm \sqrt{5}}{2}$}
\end{picture}
\begin{picture}(40,10)
\put(0,0){\circle*{3}}
\put(0,0){\line(1,0){10}}
\put(10,0){\circle*{3}}
\put(0,0){\lloop}
\put(10,0){\ulloop}
\put(10,0){\line(0,1){10}}
\put(10,10){\circle*{3}}
\put(10,0){\line(1,0){10}}
\put(20,0){\circle*{3}}
\put(-10,-10){\tiny $\lambda^3-2\lambda^2-2\lambda+2$}
\end{picture}

\bigskip
\bigskip

\begin{picture}(70,10)
\put(0,0){\circle*{3}}
\put(0,0){\line(1,0){10}}
\put(10,0){\circle*{3}}
\put(0,0){\uloop}
\put(20,0){\uloop}
\put(10,0){\line(0,1){10}}
\put(10,10){\circle*{3}}
\put(10,0){\line(1,0){10}}
\put(20,0){\circle*{3}}
\put(-20,-10){\tiny $\lambda^4-2\lambda^3-2\lambda^2+4\lambda -1$}
\end{picture}
\begin{picture}(40,10)
\put(0,0){\circle*{3}}
\put(0,0){\line(1,0){10}}
\put(10,0){\circle*{3}}
\put(0,0){\uloop}
\put(10,0){\line(0,1){10}}
\put(10,10){\circle*{3}}
\put(10,0){\line(1,0){10}}
\put(20,0){\circle*{3}}
\put(10,-10){$2$}
\end{picture}
\begin{picture}(40,10)
\put(0,0){\circle*{3}}
\put(0,0){\line(1,0){10}}
\put(10,0){\circle*{3}}
\put(10,0){\ulloop}
\put(10,0){\line(0,1){10}}
\put(10,10){\circle*{3}}
\put(10,0){\line(1,0){10}}
\put(20,0){\circle*{3}}
\put(0,-12){\small $\frac{1\pm \sqrt{13}}{2}$}
\end{picture}
\begin{picture}(40,10)
\put(0,0){\circle*{3}}
\put(0,0){\line(1,0){10}}
\put(10,0){\circle*{3}}
\put(0,0){\lloop}
\put(10,0){\ulloop}
\put(20,0){\uloop}
\put(10,0){\line(0,1){10}}
\put(10,10){\circle*{3}}
\put(10,0){\line(1,0){10}}
\put(20,0){\circle*{3}}
\put(0,-12){\small $\frac{3\pm \sqrt{5}}{2}$}
\end{picture}
\begin{picture}(40,10)
\put(0,0){\circle*{3}}
\put(0,0){\line(1,0){10}}
\put(10,0){\circle*{3}}
\put(0,0){\lloop}
\put(10,10){\lloop}
\put(20,0){\uloop}
\put(10,0){\line(0,1){10}}
\put(10,10){\circle*{3}}
\put(10,0){\line(1,0){10}}
\put(20,0){\circle*{3}}
\put(0,-12){\small $\frac{1\pm \sqrt{13}}{2}$}
\end{picture}
\begin{picture}(40,10)
\put(0,0){\circle*{3}}
\put(0,0){\line(1,0){10}}
\put(10,0){\circle*{3}}
\put(0,0){\lloop}
\put(10,10){\rloop}
\put(20,0){\rloop}
\put(10,0){\ulloop}
\put(10,0){\line(0,1){10}}
\put(10,10){\circle*{3}}
\put(10,0){\line(1,0){10}}
\put(20,0){\circle*{3}}
\put(0,-12){\small $1\pm \sqrt{3}$}
\end{picture}

\bigskip

\subsection{The roots of unity case} In this case we recover the 
Ocneanu-Kirillov-Ostrik classification of the module categories over $\cC_q$
(the quantum ``McKay correspondence''), see \cite{Oc}, \cite{KO}, \cite{O}. 
Let $q$ be a root of unity of order $N\ge 3$. Recall here that the irreducible
based modules over the based ring $\Gr (\cC_q)$ are classified by the ADET
Dynkin diagrams with the Coxeter number $N_*$, see \cite{DZ}, \cite{EK}. In the
pictures below the subscript is the number of vertices and $h$ is the
Coxeter number:
\smallskip

\noindent
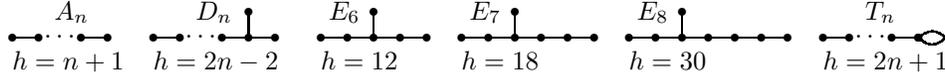
\begin{picture}(50,16)
\put(5,0){\circle*{3}}
\put(5,0){\line(1,0){10}}
\put(15,0){\circle*{3}}
\put(17,0){$\ldots$}
\put(31,0){\circle*{3}}
\put(31,0){\line(1,0){10}}
\put(41,0){\circle*{3}}
\put(21,7){$A_n$}
\put(5,-12){$h=n+1$}
\end{picture}
\begin{picture}(60,16)
\put(5,0){\circle*{3}}
\put(5,0){\line(1,0){10}}
\put(15,0){\circle*{3}}
\put(17,0){$\ldots$}
\put(31,0){\circle*{3}}
\put(31,0){\line(1,0){10}}
\put(41,0){\circle*{3}}
\put(41,0){\line(1,0){10}}
\put(41,0){\line(0,1){10}}
\put(51,0){\circle*{3}}
\put(41,10){\circle*{3}}
\put(21,7){$D_n$}
\put(5,-12){$h=2n-2$}
\end{picture}
\begin{picture}(50,16)
\put(5,0){\circle*{3}}
\put(5,0){\line(1,0){10}}
\put(15,0){\circle*{3}}
\put(15,0){\line(1,0){10}}
\put(25,0){\circle*{3}}
\put(25,0){\line(1,0){10}}
\put(25,0){\line(0,1){10}}
\put(25,10){\circle*{3}}
\put(35,0){\circle*{3}}
\put(35,0){\line(1,0){10}}
\put(45,0){\circle*{3}}
\put(8,7){$E_6$}
\put(5,-12){$h=12$}
\end{picture}
\begin{picture}(60,16)
\put(5,0){\circle*{3}}
\put(5,0){\line(1,0){10}}
\put(15,0){\circle*{3}}
\put(15,0){\line(1,0){10}}
\put(25,0){\circle*{3}}
\put(25,0){\line(1,0){10}}
\put(25,0){\line(0,1){10}}
\put(25,10){\circle*{3}}
\put(35,0){\circle*{3}}
\put(35,0){\line(1,0){10}}
\put(45,0){\circle*{3}}
\put(45,0){\line(1,0){10}}
\put(55,0){\circle*{3}}
\put(8,7){$E_7$}
\put(5,-12){$h=18$}
\end{picture}
\begin{picture}(70,16)
\put(5,0){\circle*{3}}
\put(5,0){\line(1,0){10}}
\put(15,0){\circle*{3}}
\put(15,0){\line(1,0){10}}
\put(25,0){\circle*{3}}
\put(25,0){\line(1,0){10}}
\put(25,0){\line(0,1){10}}
\put(25,10){\circle*{3}}
\put(35,0){\circle*{3}}
\put(35,0){\line(1,0){10}}
\put(45,0){\circle*{3}}
\put(45,0){\line(1,0){10}}
\put(55,0){\circle*{3}}
\put(55,0){\line(1,0){10}}
\put(65,0){\circle*{3}}
\put(8,7){$E_8$}
\put(5,-12){$h=30$}
\end{picture}
\begin{picture}(50,12)
\put(5,0){\circle*{3}}
\put(5,0){\line(1,0){10}}
\put(15,0){\circle*{3}}
\put(17,0){$\ldots$}
\put(31,0){\circle*{3}}
\put(31,0){\line(1,0){10}}
\put(41,0){\circle*{3}}
\put(41,0){\begin{picture}(10,10)
\put(0,0){\circle*{3}}
\qbezier(0,0)(5,5)(10,0)
\qbezier(0,0)(5,-5)(10,0)
\end{picture}}
\put(21,7){$T_n$}
\put(5,-12){$h=2n+1$}
\end{picture}
\bigskip

\begin{theorem} Let $q$ be a primitive root of unity of order $N\ge 3$. 

{\em (i)} Assume that $N$ is even. The indecomposable module categories over 
the category $\cC_q$ are classified by the ADE Dynkin diagrams with the 
Coxeter number $N_*=N/2$.

{\em (ii)} Assume that $N$ is odd. The indecomposable module categories over 
the category $\cC_q$ are classified by the ADET Dynkin diagrams with the 
Coxeter number $N_*=N$.
\end{theorem}

\begin{proof} It is clear that any module category $\cM$ over $\cC_q$ gives 
rise to a based module $\Gr (\cM)$ over $\Gr (\cC_q)$. Such based modules
were classified in \cite{DZ}, \cite{EK} and the answer is given precisely by 
ADET Dynkin diagrams with the Coxeter number $N_*$. Conversely, we know that 
any generalized tree with nondegenerate eigenvalue $q+q^{-1}$ gives rise to a
unique module category over $\cT_q$. It is well known that the class of the
object $\Ima (f_{N_*-1})$ in $\Gr (\tilde \cC_q)$ is given by $P_{N_*}([V])$,
see e.g. \cite{BK} (recall that $P_{N_*}$ is an ultraspherical polynomial). 
On the other hand in $\Gr (\cC_q)$ we have the relation 
$P_{N_*}([V])=0$, see {\em loc. cit.} Thus we can apply Theorem \ref{turaev2}.
Observe that $-q-q^{-1}$ is a nondegenerate eigenvalue of the adjacency
matrix of the corresponding graph $\Gamma$ in all cases except when $N$ is 
even and $\Gamma =T_n$ (actually in all cases $-q-q^{-1}$ is Galois conjugate
to the Frobenius-Perron eigenvalue). Thus by Proposition \ref{nond} we have
a unique solution of the main equation. The Theorem is proved.
\end{proof}


\begin{thebibliography}{999}

\bibitem[BK]{BK} B.~Bakalov, A.~Kirillov Jr.,
\textit{Lectures on Tensor categories and modular functors}, 
AMS, Providence, (2001).

\bibitem[Bi]{Bi} J.~Bichon, {\em The representation category of the quantum
group of a non-degenerate bilinear form}, math.QA/0111114.

\bibitem[B]{B} A.~Bondal, {\em A symplectic groupoid of triangular bilinear
forms and the braid group}, preprint IHES, available at
www.ihes.fr/PREPRINTS/M00/Resu/resu-M00-02.html

\bibitem[DZ]{DZ} P.~Di~Francesco, J.-B.~Zuber, {$SU(N)$ lattice integrable 
models associated with graphs}, Nucl. Phys. B 338 (1990), 602-646.

\bibitem[D]{D} J.~Du, {\em A note on quantized Weyl reciprocity at roots
of unity}, Algebra Colloq. {\bf 2} (1995), no. 4, 363-372.

\bibitem[EK]{EK} P.~Etingof, M.~Khovanov, {\em Representations of tensor 
categories and Dynkin diagrams},. Internat. Math. Res. Notices 1995, no. 5, 
235--247.

\bibitem[F]{Fi} M.~Fiedler, {\em Eigenvectors of acyclic matrices}, 
Czechoslovak Mathematical Journal {\bf 25} (1975), 607-618.

\bibitem[K]{K} C.~Kassel, {\em Quantum groups}, Graduate Texts in Mathematics
{\bf 155}, Springer-Verlag, New York, 1995.

\bibitem[KO]{KO} A.~Kirillov Jr., V.~Ostrik, {\em On a $q$-analogue of the 
McKay correspondence and the ADE classification of 
$sl_2$ conformal field theories}, Adv. Math.
{\bf  171} (2002), no. 2, 183--227.

\bibitem[Oc]{Oc} A.~Ocneanu, {\em The classification of subgroups of quantum
$SU(N)$}, in Quantum symmetries in theoretical physics and mathematics
(Bariloche, 2000), Contemp. Math. {\bf 294} (2002), 133-159.

\bibitem[O]{O} V.~Ostrik, {\em Module categories, weak Hopf algebras and 
modular invariants}, preprint math.QA/0111139.

\bibitem[T]{T} V.~Turaev, {\em Quantum invariants of knots and 3-manifolds},
de Gruyter Studies in Mathematics, 18. Walter de Gruyter \& Co., Berlin, 1994.
\end{thebibliography}
\end{document}